\documentclass[12pt]{article}
\usepackage{latexsym}
\usepackage{amsmath}
\usepackage{amssymb}
\usepackage{amsfonts}
\usepackage{color}
\usepackage{multirow}
\usepackage{srcltx}
\def\norm2to2{{\|\cdot\|_{2,2}}}
\def\Prob{\hbox{\rm Prob}}

\def\bE{{\mathbf{E}}}

\def\S{{\mathbf{S}}}
\def\Diag{\hbox{\rm  Diag}}
\def\Prob{\hbox{\rm  Prob}}

\def\Opt{\hbox{\rm Opt}}

\def\Tr{{\mathop{\hbox{\rm  Tr}}}}
\def\cA{{\cal A}}
\def\cB{{\cal B}}

\def\cE{{\cal E}}

\def\cN{{\cal N}}

\def\cS{{\cal S}}

\def\Argmin{\mathop{\hbox{\rm  Argmin}}}

\def\Det{{\mathop{\hbox{\rm  Det}}}}

\def\S{{\mathbf{S}}}

\def\R{{\mathbb R}}
\def\C{{\mathbb C}}
\def\Z{{\mathbb Z}}

\def\qed{\ \hfill$\square$\par\smallskip}
\newcommand{\be}{\begin{eqnarray}}
\newcommand{\ee}[1]{\label{#1}\end{eqnarray}}

\newcommand{\ese}{\end{eqnarray*}}
\newcommand{\bse}{\begin{eqnarray*}}
\newcommand{\rf}[1]{~(\ref{#1})}

\def\mypict3{\epsfxsize=220pt\epsfysize=80pt\epsfbox}

\newtheorem{lemma}{Lemma}[section]
\newtheorem{proposition}{Proposition}[section]

\newtheorem{example}{Example}[section]

\title{On Low Rank Matrix Approximations with Applications to Synthesis Problem in Compressed Sensing}
\author{ Anatoli Juditsky\thanks{LJK,
Universit\'e J. Fourier, B.P. 53, 38041 Grenoble Cedex 9, France, {\tt Anatoli.Juditsky@imag.fr}} \and \and Fatma K{\i}l{\i}n\c{c} Karzan\thanks{Georgia Institute
 of Technology, Atlanta, Georgia
30332, USA, {\tt fkilinc@isye.gatech.edu}} \and Arkadi Nemirovski\thanks{Georgia Institute
 of Technology, Atlanta, Georgia
30332, USA, {\tt{nemirovs@isye.gatech.edu}}. Research of the second and the third authors was supported by the Office of Naval Research grant \# N000140811104.}}
\begin{document}
\maketitle
\begin{abstract}
We consider the {\sl synthesis problem} of Compressed Sensing --given $s$ and an $M\times n$ matrix $A$,  extract from it an $m\times n$  submatrix $A_m$, certified to be $s$-good, with $m$ as small as possible. Starting from the verifiable sufficient conditions of $s$-goodness, we express the synthesis problem as the problem of approximating a given matrix by a matrix of specified low rank in the uniform norm. We propose randomized algorithms for efficient construction of rank $k$ approximation of matrices of size $m\times n$ achieving accuracy bounds $O(1)\sqrt{{\ln(mn)\over k}}$ which hold in expectation or with high probability. We also supply derandomized versions of the approximation algorithms which does not require random sampling of matrices and attains the same accuracy bounds. We further demonstrate that our algorithms are optimal up to the logarithmic in $m,n$ factor, i.e. the accuracy of such an approximation for the identity matrix $I_n$ cannot be better than $O(1)k^{-{1\over 2}}$. We provide preliminary numerical results on the performance of our algorithms for the synthesis problem.
\end{abstract}
\section{Introduction}
Let $A\in\R^{m\times n}$ be a matrix with $m<n$. Compressed Sensing focuses on recovery of a {\sl sparse} signal $x\in\R^n$ from its noisy observations
$$
y=Ax+\delta,
$$
where $\delta$ is an observation noise such that $\|\delta\|\leq\epsilon$ for certain known norm on $\R^m$ and some given $\epsilon$. The standard recovering routine is
$$
\widehat{x}\in\Argmin_w\{\|w\|_1: \|Aw-y\|\leq\epsilon.\}.
$$
We call the matrix $A$ $s$-{\sl good} if whenever the true signal $x$ is $s$-sparse (i.e., has at most $s$ nonzero entries) and there is no observation errors  ($\epsilon=0$),  $x$ is the unique optimal solution to the optimization program $\min\{\|w\|_1:Aw=Ax\}.$
\par
To the best of our knowledge, nearly the strongest verifiable sufficient condition for $A$ to be $s$-good is as follows (cf \cite{CSFirst}): \be {\rm There \;exists}\; Y\in\R^{m\times n} \;{\rm such \;that} \;\|I_n-Y^TA\|_\infty<{1\over 2s} \ee{good1} (here and in what follows $\|X\|_\infty=\max\limits_{i,j}|X_{ij}|$, $X_{ij}$ being the elements of $X$).\footnote{ We address the reader to \cite{CSFirst} for details concerning the derivation, the link to the necessary and sufficient condition of $s$-goodness and its comparison to traditional non-verifiable sufficient conditions for
$s$-goodness based on Restricted Isometry or Restricted Eigenvalue Property and a verifiable sufficient condition based on mutual incoherence.}
\par
In this paper we consider the {\sl synthesis problem} of Compressed Sensing as follows:
\begin{quotation}
{\em Given $s$ and an  $M\times n$ matrix $A$,  extract from it an $m\times n$  submatrix $A_m$, certified to be $s$-good, with $m$ as small as possible.}
\end{quotation}
One can think, e.g., of a spatial or planar $n$-point grid $\cE$ of possible locations of signal sources and an $M$-element grid $\cS$ of possible locations of sensors. A sensor in a given location measures a known, depending on the location, linear form of the signals emitted at the nodes of $\cE$, and the goal is to place a given number $m \ll M$ of sensors at the nodes of $\cS$ in order to be able to recover the location of sources via the $\ell_1$-minimization, conditioned that there are $s$ sources at most.
Since the property of $s$-goodness is difficult to verify, we will look for a submatrix of the original matrix $A$ for which the $s$-goodness can be certified by the sufficient condition \rf{good1}.
Suppose that along with $A$ we know an $M\times n$ matrix $Y_M$ which certifies that the ``level of goodness'' of $A$
is at least $s$, that is, we have
\begin{equation}\label{invok}
\|I_n-\bar Y_M^TA\|_\infty\le \mu <{1\over 2s}.
\end{equation}
Then we can approach the synthesis problem as follows:
\begin{quotation}
Given $M\times n$ matrices $Y_M$ and $A$ and a tolerance $\epsilon>0$, we want to extract from $A$, $m$ rows (the smaller is $m$, the better) to get an $m\times n$ matrix $A_m$ which, along with properly chosen $Y_m\in\R^{m\times n}$, satisfies the relation $\|Y_M^TA-Y_m^TA_m\|_\infty\leq\epsilon$.
\end{quotation}
Choosing $\epsilon<{1\over 2s}-\mu$ and invoking (\ref{invok}), we ensure that the output $A_m$ of the above procedure is $s$-good. This simple observation motivates our interest to the problem of approximating a given matrix by a matrix of specified (low rank) in the {\em uniform norm}.

Note that in the existing literature on low rank approximation of matrices the emphasis is on efficient construction when the approximation error is measured in the Frobenius norm (for the Frobenius norm $\|A\|_F=\left(\sum_{i,j}A_{ij}^2\right)^{1/2}$). Though the Singular Value Decomposition (SVD) gives the best rank $k$ approximation in terms of all the norms that are invariant under rotation (e.g., the Frobenius norm and the spectral norm), its computational cost may be prohibitive for applications involving large matrices. Recently, the properties of fast low rank approximations in the Frobenius norm based on the randomized sampling of rows (or columns) of the matrix  (see, e.g., \cite{Drineas_Kannan_Mahoney_06_ii,Frieze_Kannan_Vempala_04}) or random sampling of a few individual entries (see \cite{Achlioptas_McSherry_07} and references therein) has been studied extensively. Another randomized fast approximation based on the preprocessing by the Fast Fourier Transform or Fast Hadamard Transform has been studied in \cite{Nguyen_Do_Tran_09}. Yet we do not know explicit bounds available from the previous literature which concern numerically efficient low rank approximations in the uniform norm.

In this work, we aim at developing efficient algorithms for building low rank approximation of a given matrix in the uniform norm. Specifically, we consider two types of low rank approximations:
\begin{enumerate}
\item Let $W=Y^TA$, where $Y$ and $A$ are known $M\times n$ matrices. We consider the approximation $W_k=Y_kA^T_k$ of $W$ such that the matrices $Y_k$ and $A_k$ of dimension $m_k\times n$, $m_k\le k\le M$, are composed of  multiples of the rows of the matrices $Y$ and $A$ respectively. We show that a fast (essentially, of numerical complexity $O(kMn^2)$) approximation $W_k$ can be constructed which satisfies
\[
\|W-W_k\|_\infty=O(1)L(Y,A)\sqrt{\ln [n]\over k},
\]
where $L(Y,A)=\sum_{i}\|y_i\|_\infty\|a_i\|_\infty$ and $y_i^T,a_i^T$ denote the $i$-th rows of $Y$ and $A$ respectively. Note that for moderate values of $L(Y,A)=O(1)$ and $k<n/2$ this approximation is ``quasi-optimal'', as we know (cf, e.g. \cite[Proposition 4.2]{CSFirst}) that (for certain matrices $W$) the accuracy of such an approximation cannot be better than $O(k^{-1/2})$. \item Let $A\in \R^{m\times n}$, $A=MN^T$, where $M\in
\R^{m\times d}$ and $N\in\R^{n\times d}$. We consider a fast  approximation $A_k=\sum_{i=1}^k \eta_i\zeta_i^T$ of $A$, where $\eta_i$ and $\zeta_i$ are linear combinations of columns of $M$ and $N$ respectively. We show that this approximation satisfies
\[
\|A-A_k\|_\infty\le O(1)D\sqrt{\ln [mn]\over k}
\]
where $D$ is the maximal Euclidean norm of rows of $M$ and $N$. We show that when $A$ is an $n\times n$  identity matrix  the above bound is unimprovable up to a logarithmic factor.
\end{enumerate}
In this paper we propose two types of construction of fast approximations: we consider the randomized construction, for which the accuracy bounds above hold in expectation (or with significant
probability). We also supply ``derandomized'' versions of the approximation algorithms which does not require random sampling of matrices and attains the same accuracy bounds as the randomized method.

\section{Low rank approximation in Compressed Sensing}
In this section we suppose to be given $s$ and an  $M\times n$ matrix $A$ and our objective  is to extract from $A$ a submatrix $A_k$ which is composed of, at most, $k$ rows of $A$, with as small $k$ as possible, which is $s$-good. We assume that $A$ admits a ``goodness certificate'' $Y$. Namely, we are given an $M\times n$ matrix $Y$ such that \be \mu:=\|I_n-Y^TA\|_\infty< {1\over 2s}, \ee{eq20} and we are looking for $A_k$ and the corresponding $Y_k$ such that $\|I_n-Y_k^TA_k\|< {1\over 2s}$.
\subsection{Random sampling algorithm}\label{RandSamp}
The starting point of our developments is the following simple
\begin{lemma}\label{dentropy}
Let for $\beta >0$, let \be V_\beta(z)=\beta \ln \left[\sum_{i=1}^d\cosh\left({z_i\over \beta}\right)\right]-\beta\ln d:\R^{d}\times \R_+\to \R_+. \ee{Vb}
Then
\begin{itemize}
\item [(i)] we have $\|z\|_\infty-\beta\ln [2d]\le V_\beta(z)\le \|z\|_\infty$; \item[(ii)] if $\beta_1\le \beta_2$ then $V_{\beta_1}(z)\ge V_{\beta_2}(z)$; \item[(iii)] function $V_\beta$ is convex and continuously differentiable on $\R^{d}$. Further, its gradient $V'_\beta$ is Lipschitz-continuous with the constant $\beta^{-1}$: \be \|V'_\beta(z_1)-V'_\beta(z_2)\|_1\le\beta^{-1}\|z_1-z_2\|_\infty, \ee{2td} and $\|V'_\beta(z)\|_1\le 1$ for all $z\in \R^{d}$.
 \end{itemize}
\end{lemma}
For proof, see Appendix \ref{appa}.

Lemma \ref{dentropy} has the following immediate consequence:
\begin{proposition}\label{basic}
Let $\beta\ge \beta'>0$ (non-random) and let $\xi_1$,...,$\xi_k$ be random vectors in $\R^d$ such that  $\bE\{\xi_i|\xi_1,...,\xi_{i-1}\}=0$ a.s., and $\bE\{\|\xi_i\|_\infty^2\}\leq\sigma_i^2<\infty$ for all $i \in \{1,\ldots,k\}$,
and let $S_k=\sum_{i=1}^k \xi_k$. Then \be \bE\{V_\beta(S_k)\}\le \bE\{V_{\beta'}(S_{k-1})\}+{\sigma^2_k\over 2\beta'}. \ee{derand0} As a result, \be \bE\left\{\left\|S_k\right\|_\infty\right\}\leq
\sqrt{2\ln [2d]\sum_{i=1}^k\sigma_i^2}. \ee{eq7}
\end{proposition}
{\bf Proof.} Let $\beta\ge \beta'$. By applying items (ii) and  (iii) of the lemma we get:
\[
V_{\beta}(S_k)\le V_{\beta'}(S_k)\le V_{\beta'}(S_{k-1})+\langle V'_{\beta'}(S_{k-1}),\,\xi_k\rangle+{1\over 2{\beta'}}\|\xi_k\|^2_\infty.
\]
When taking the expectation (first conditional to $\xi_1,...,\xi_{k-1}$), due to $\bE \{\xi_k|\xi_1,...,\xi_{k-1}\}=0$ a.s.,  we obtain
\[
\bE \{V_{\beta}(\S_k)\}\le \bE \{V_{\beta'}(S_{k-1})\}+{\bE \{\|\xi_k\|_\infty^2\}\over 2{\beta'}}\le\bE \{V_{\beta'}(S_{k-1})\}+{\sigma_k^2\over 2{\beta'}},
\]
which is \rf{derand0}. Now let us set $\beta'=\beta=\sqrt{\sum_{i=1}^k \sigma^2_i\over 2\ln[2d]}$. Since $V_{\beta}(0)=0$ we conclude that $$\bE \{V_{\beta}(S_k)\}\le\sum_{i=1}^k{\sigma_i^2\over 2{\beta}}.$$ On the other hand, by item (i) of Lemma \ref{dentropy},
\bse \bE \{\|S_k\|_\infty\}\le {\beta}\ln [2d]+\bE \{V_{\beta}(S_k)\}\le \beta\ln [2d]+\sum_{i=1}^k{\sigma^2_i\over 2{\beta}}\le
\sqrt{2\ln [2d]\sum_{i=1}^k\sigma_i^2}. \ese
 \qed
\paragraph{The random sampling algorithm.} Denoting $y^T_i$ and $a^T_i,\,i=1,...,M$, $i$-th rows of $Y$ and  $A$, respectively,  let us set
\begin{equation}\label{L}
\theta_i=\|y_i\|_\infty\,\|a_i\|_\infty, \;\; L=\sum_i\theta_i,\;\;\pi_i={\theta_i\over L}, \,\,z_i={L\over \theta_i}y_i,
\end{equation}
and let $W=Y^TA$. Observe that \be
\begin{array}{rcl}
W&=&\sum_{i=1}^M \pi_i\left(z_ia_i^T\right),\\ \|z_ia_i^T\|_\infty&=&L,\;\;1\leq i\leq M,\\ \sum_{i=1}^M\pi_i&=&1,\;\;\pi_i\geq 0,\;\;1\leq i\leq M.
\end{array}
\ee{then} Now let $\Xi$ be random rank 1 matrix taking values $z_ia_i^T$ with probabilities $\pi_i$, and let $\Xi_1,\Xi_2,...$ be a sample of independent realizations of $\Xi$. Consider the random matrix
$$
W_k={1\over k} \sum_{\ell=1}^k\Xi_\ell.
$$
Then $W_k$ is, by construction, of the form $Y_k^TA_k$, where $A_k$ is  a random $m_k\times n$ submatrix of $A$ with $m_k\leq k$.

As an immediate consequence of Proposition \ref{basic} we obtain the following statement:
\begin{proposition}\label{obs30} One has
\begin{equation}\label{eq40}
\bE\left\{\|W_k-W\|_\infty\right\}\leq 2Lk^{-1/2}\sqrt{2\ln(2n^2)}.
\end{equation}
In particular, the probability of the event
$$
\cE=\{\Xi_1,...,\Xi_k:\|W_k-W\|_\infty\leq 4Lk^{-1/2}\sqrt{2\ln[2n^2]}\}
$$
is $\geq1/2$, and whenever this event takes place, we have in our disposal a matrix $Y_k$ and a $m_k\times n$ submatrix $A_k$ of $A$ with $m_k\leq k$ such that
\begin{equation}\label{suchthat}
\|I_n-Y_k^TA_k\|_\infty\leq \|I_n-W\|_\infty+\|W_k-W\|_\infty\leq \mu_k:=\mu+ 4Lk^{-1/2}\sqrt{2\ln[2n^2]}.
\end{equation}
\end{proposition}
{\bf Proof.} By (\ref{then}) we have $\|z_ia_i^T\|_\infty\leq L$ for all $i$, and besides this, treating $i$ as random index distributed in $\{1,...,n\}$ according to probability distribution $\pi=\{\pi_i\}_{i=1}^n$, we have $\bE\{z_ia_i^T\}=W$. It follows that $\|\Xi_\ell-W\|_\infty\leq 2L$ and $\bE\{\Xi_\ell-W\}=0$. If we  denote $S_i=\sum_{\ell=1}^i (\Xi_\ell-W)$, when applying Lemma \ref{dentropy}  we obtain \bse \bE \{\|S_k\|_\infty\}\le 2L\sqrt{2k\ln [2n^2]}, \ese and we arrive at \rf{eq40}.
\qed
\paragraph{Discussion.} Proposition \ref{obs30} suggests a certain approach to the synthesis problem. Indeed, according to this Proposition, picking at random $k$ rows $a^T_{i_\ell}$, where $i_1,...,i_k$ are sampled independently from the distribution $\pi$, we get with probability at least $1/2$ a random $m_k\times n$ matrix $A_k$, $m_k\leq k$, which is provably $s$-good with $s=O(1)(L\sqrt{\ln[n]/k}+\mu)^{-1}$. When $L=O(1)$, this is nearly as good as it could be, since the sufficient condition for $s$-goodness stated in \rf{good1} can justify $s$-goodness of an $m\times n$ sensing matrix with $n>O(1)m$ only when $s\leq O(1)\sqrt{m}$, see \cite[Proposition 4.2]{CSFirst}. %
\subsection{Derandomization}
\label{derandom} Looking at the proof of Proposition \ref{basic}, we see that the construction of $A_k$ and $Y_k$ can be derandomized. Indeed, \rf{derand0} implies that
\begin{quote}
{\sl Whenever $S\in\R^{n\times n}$ and $\beta\ge \beta'$ there exists $i$ such that
\[
V_{\beta}(S+(z_ia_i^T-W))\leq V_{\beta'}(S)+{2L^2\over \beta'}.
\]
Specifically, the above bound is satisfied for every $i$ such that}
\[
\langle V'_{\beta'}(S),z_ia_i^T-W\rangle\leq 0,
\]
and because $\pi_i \ge 0$ and $\sum_i\pi_i(z_ia_i^T-W)=0$, the latter inequality is certainly satisfied for some $i$.
\end{quote}
Now assume that given a sequence $\beta_0\ge\beta_1\ge ...$ of positive reals, we build a sequence of matrices $S_i$ according to the following rules:
\begin{enumerate}
\item $S_0=0$; \item $S_{k+1}=S_k+(v_ka_{\ell_k}^T-W)$ with $\ell_k\in\{1,...,M\}$ and $v_k\in\R^n$ such that
\begin{equation}\label{condition}
V_{\beta_{k+1}}(S_{k+1})\leq V_{\beta_k}(S_k)+\delta_k,\;\;\delta_k\leq {2L^2\over \beta_k}.
\end{equation}
\end{enumerate}
Then for every $k\geq1$ the matrix $U_k=k^{-1}S_k$ is of the form $Y_k^TA_k-W$, where $A_k$ is a $m_k\times n$ submatrix of $A$ with $m_k\leq k$, and
$$
\|S_k\|_\infty\leq \beta_{k}\ln[2n^2]+\sum_{\ell=0}^{k-1}\delta_\ell,
$$
whence
\[
\|Y_k^TA_k-I_n\|_\infty\leq \mu+k^{-1}\left(\beta_{k}\ln[2n^2]+\sum_{\ell=1}^k\delta_\ell \right).
\]
In particular, for the choice $\beta_\ell= L\sqrt{2k\over \ln [2n^2]}$, $\ell=0,...,k$, we obtain
\[
\|Y_k^TA_k-I_n\|_\infty\leq\mu+2L\sqrt{2\ln [2n^2]\over k}
\]
One can consider at least the following  three (numerically efficient) policies for choosing $v_k$ and $\ell_k$ satisfying (\ref{condition}); we order them according to their computational complexity.
\begin{description}
\item[A.] Given $S_k$, we test one by one the options $\ell_k=i$, $v_k=z_i$, $i=1,...,M$, until an option satisfying (\ref{condition}) is met (or test all the $n$ options and choose the one which results in the smallest $V_{\beta_{k+1}}(S_{k+1})$). Note that accomplishing a step of this scheme requires $O(Mn^2)$ elementary operations.
\par
\item [A$'$.] In this version of A, we test the options $\ell_k=i$, $v_k=z_i$ when picking $i$ at random, as independent realizations of the random variable $\imath$ taking values $1,...,M$ with probabilities $\pi_i$, until an option with $\langle V'_{\beta_k}(S_k),z_ia_i^T-W\rangle\leq0$ is met. Since $\bE\left\{\langle V'_{\beta_k}(S_k),z_ia_i^T-W\rangle\right\}\le 0$, we may  hope that this procedure will take essentially less steps than the ordered scan through the entire range $1,...,M$ of values of $i$. \item [B.] Given $S_k$ we solve $M$ one-dimensional convex optimization problems \be t_i^*\in\Argmin_{t\in \R_+} V_{\beta_{k}}(S_k+tz_ia_i^T-W),\,\,1\leq i\leq M, \ee{bisec} then select the one, let its index be $i_*$,  with the smallest value of $
    V_{\beta_{k}}(S_k+t_i^*z_ia_i^T-W)$, and put $v_k=t_{i_*}^*z_{i_*}$, $\ell_k=i_*$.

If the bisection algorithm is used to find $t^*_i$, solving the problem \rf{bisec} for one $i$ to the relative accuracy $\epsilon$ requires $O(n^2\ln \epsilon^{-1})$ elementary operations. The total numerical complexity of the step of the method is $O(Mn^2\ln \epsilon^{-1})$. \item [C.] Given $S_k$, we solve $M$ convex optimization problems \be u_i^*\in\Argmin_{u\in
\R^n}V_{\beta_k}(S_k+ua_i^T-W),\;\;1\le i\le M, \ee{bisec2} then select the one, let its index be $i_*$,  with the smallest value of $V_{\beta_k}(S_k+u_i^*a_i^T-W)$, and set $v_k=u_i^*$, $\ell_k=i_*$.

Note that due to the structure of $V_\beta$  to solve \rf{bisec2} it suffices to find a solution to the system
\be
\begin{array}{l}\sum_{\ell=1}^n \gamma_\ell\sinh(\alpha_{j\ell}+\gamma_{\ell}u_j)=0,\\ \alpha_{j\ell}={[S_k]_{j\ell}-[W]_{j\ell}\over \beta_k},\;\;\gamma_{\ell}={[A]_{\ell i}\over \beta_k}, \;\;1\le j,\ell\le n.\end{array} \ee{bisec3}
Since the equations of the system \rf{bisec3} are independent, one can use bisection to find the component $u_j$ of the solution.\footnote{Note that due to the convexity of the left-hand side of the equation in \rf{bisec3}, even faster algorithm of Newton family can be used.} Finding a solution to the relative accuracy $\epsilon$ to each equation then requires $O(n\ln \epsilon^{-1})$ arithmetical operations, and the total complexity of solving \rf{bisec2} becomes $O(Mn^2\ln \epsilon^{-1})$.
\end{description}
\paragraph{Selecting  $Y$ and $W$.} Note that the numerical schemes of this section should be initialized with matrices $Y$ and $W=Y^TA$. We can do as follows:
\begin{enumerate}
\item We start with solving the problem
$$
Y\in\Argmin_{Z=[z_1^T;...;z_M^T]\in\R^{M\times n}}\left\{\sum_{i=1}^M\|z_i\|_\infty\|a_i^T\|_\infty:\; \|I_n-Z^TA\|_\infty\leq\mu\right\},
$$
where $\mu$ is a certain fraction of ${1\over 2s}$. Assuming the problem is feasible for the chosen $\mu$, we get in this way the ``initial point'' -- the matrix $W=Y^TA$. \item Then we apply the outlined procedure to find $A_k$ and  $Y_k$. At each step $\ell$ of this procedure, we get certain $m_\ell\times n$ submatrix $A_\ell$ of $A$ and a matrix $Y_\ell$. When
$\|I_n-Y_\ell^TA_\ell\|_\infty$ becomes less than ${1\over 2s}$ we terminate. Alternatively, we can solve at each step $\ell$ an auxiliary problem $\min\limits_{U\in \R^{m_\ell\times n}}\|I_n-U^TA_\ell\|_\infty$ and terminate when the optimal value in this problem becomes less than ${1\over 2s}$.
\end{enumerate}
\paragraph{Choosing the sequence $(\beta_\ell)$.} When the number $k$ of steps of the iterative schemes of this section is fixed, the proof of Proposition \ref{basic} suggests the fixed choice of the ``gain sequence'' $(\beta_\ell)$: $\beta_\ell=L\sqrt{2k\over \ln[2n^2]}$, $\ell=1,...,k$. When the number $k$ is not known {\em a priori}, one can use the sequence, computed recursively according to the rule $\beta_{\ell}=\beta_{\ell-1}+{2L^2\over \ln[2n^2]\beta_{\ell-1}}$, $\beta_0={2L^2\over \ln[2n^2]}$, or, what is essentially the same,  the sequence $\beta_\ell=2L\sqrt{\ell+1\over \ln [2n^2]}$, $\ell=0,1,...$.  Another possible choice of $\beta_\ell$'s is as follows: observe first that the function $V_\beta(z)$ is jointly convex in $\beta$ and $z$. Therefore, we may modify the above algorithms by adding the minimization in $\beta$. For instance, instead of the optimization problems \rf{bisec} in item B we can consider $M$ two-dimensional optimization problems
\bse (t_i^*,\beta_i^*)\in\Argmin_{t,\beta\in \R_+}\left\{ \beta \ln[2n^2]+V_{\beta}(S_k+tz_i[A^T]_i^T-W)\right\},\,\,1\leq i\leq M; \ese
we select the one with the smallest value of the objective $ V_{\beta_i^*}(S_k+t_i^*z_ia_i^T-W)+\beta_i^* \ln[2n^2]$, and set, as before, $v_k=t_{i_*}^*z_{i_*}$, $\ell_k=i_*$. Note that such a modification does not increase significantly the complexity estimate of the scheme.
\subsection{Numerical illustration}
Here we report on preliminary numerical experiments with the synthesis problem as posed in the introduction. In our experiment, $A$ is square, specifically, this is the Hadamard matrix $H_{11}$ of order 2048.
\begin{quote}
Recall that the Hadamard matrix $H_\nu$, $\nu=0,1,...$ is a square matrix of order $2^\nu$ given by the recurrence
$$
H_0=1,\,H_{s+1}=\left[\begin{array}{rr}H_s&H_s\cr
H_s&-H_s\cr\end{array}\right],
$$
whence $H_\nu$ is a symmetric matrix with entries $\pm1$ and $H_\nu^TH_\nu=2^\nu I_{2^\nu}$.
\end{quote} The goal of the experiment was to extract from $A=H_{11}$ an $m\times 2048$ submatrix $A_m$ which satisfies the relation (cf. (\ref{good1}))
\begin{equation}\label{relation}
\Opt(A_m):=\min\limits_{Y_m\in\R^{m\times n}} \|I_n-Y_m^TA_m\|_\infty<{1\over 2s},\,\,n=2048
\end{equation}
with $s=10$; under this requirement, we would like to have $m$ as small as possible. In Compressed Sensing terms, we are trying to solve the synthesis problem with  $A=H_{11}$; in low rank approximation terms, we want to approximate $I_{2048}$ in the uniform norm within accuracy $<0.05$ by a rank $m$ matrix of the form $Y_m^TA_m$, with the rows of $A_m$ extracted from $H_{11}$. The
advantages of the Hadamard matrix in our context is twofold:
\begin{enumerate}
\item The error bound (\ref{eq40}) is proportional to the quantity $L$ defined in (\ref{L}). By the origin of this quantity, we clearly have $\|Y^TA\|_\infty\leq L$, whence $L\geq1-\mu>1-{1\over 2s}\geq 1/2$ by (\ref{eq20}). On the other hand, with $A=H_\nu$ being an Hadamard matrix, setting $Y=2^{-n}YH_\nu$, so that $Y^TA-I_{2^\nu}$, we ensure the validity of (\ref{eq20}) with $\mu=0$ and get $L=1$, that is, $\mu$ is as small as it could be, and $L$ is nearly as small as it could be.  \item Whenever $A_m$ is a submatrix of $H_\nu$, the optimization problem in the left hand side of (\ref{relation}) is easy to solve.
\end{enumerate}
Item 2 deserves an explanation. Clearly, the optimization program in (\ref{relation}) reduces to the series of $n=2048$ LP programs
\begin{equation}\label{Opti}
\Opt_i(A_m)=\min_{y\in\R^m}\|e_i-A_m^Ty\|_\infty,\,1\leq i\leq n,
\end{equation}
where $e_i$ is the standard basic orth in $\R^n$; and $\Opt(A_m)=\max\limits_i\Opt_i(A_m)$. The point is (for justification, see Appendix \ref{appHad}) that {\sl when $A_m$ is an $m\times n$  submatrix of the $n\times n$ Hadamard matrix, $\Opt_i(A_m)$ is independent of $i$,} so that checking the inequality in (\ref{relation}) requires solving a {\sl single} LP program with $m$ variables rather than solving $n$ LO programs of the same size.
\par
The experiment was organized as follows. As it was already mentioned, we used $\nu=11$ (that is, $n=2048$) and $s=10$  (that is, the desired uniform norm of approximating $I_{2048}$ by $Y_m^TA_m$ was 0.05). We compared two approximation policies:
\begin{itemize}
\item ``Blind'' approximation -- we choose a random permutation $\sigma(\cdot)$ of the indices $1,...,2048$ and look at the submatrices $A^k$, $k=1,2,...$ obtained by extracting from $H_{11}$ rows with indices $\sigma(1),\sigma(2),...,\sigma(k)$ until a submatrix satisfying (\ref{relation}) is met. This is a refinement of the Random sampling algorithm as applied to $A=H_{11}$ and $Y=2^{-11}A$, which results in $W=I_{2048}$.   The refinement is that instead of looking for approximation of $W=I_{2048}$ of the form ${1\over k}\sum_{\ell=1}^kz_{i_\ell}a_{i_\ell}^T$, where $i_1,i_2,...$ are independent realizations of random variable $\imath$ taking values $1,...,\mu$ with equal probabilities (as prescribed by (\ref{L}) in the case of $A=H_\nu$), we look for the best approximation of the form $Y_k^TA^k$, where $A^k$ is the submatrix of $A$ with the row indices $\sigma(1),...,\sigma(k)$.
\item ``Active'' approximation, which is obtained from algorithm {\bf A$'$} by the same refinement as in the previous item.
\end{itemize}
In our experiments, we ran every policy 6 times. The results were as follows:\\
\indent
``Blind'' policy $\cB$: the rank of $0.05$-approximation of $W=I_{2048}$ varied from 662 to 680.\\
\indent ``Active'' policy $\cA$: the rank of $0.05$-approximation of $W$ varied from 617 to 630.\\
Note that in both algorithms the resulting matrix $A_m$ is built ``row by row'', and the certified levels of goodness of the intermediate matrices $A^1,A^2,...$ are computed. In the below table we indicate, for the most successful (resulting in the smallest $m$) of the 6 runs of each algorithm, the smallest values of $k$ for which $A^k$ was certified to be $s$-good, $s=1,2,...,10$:
$$
\begin{array}{||r||c|c|c|c|c|c|c|c|c|c||}
\hline\hline
s&1&2&3&4&5&6&7&8&9&10\\
\hline\hline
\cB&15&58&121&197&279&343&427&512&584&662\\
\hline
\cA&12&47&104&172&246&323&399&469&547&617\\
\hline\hline
\end{array}
$$
Finally, we remark that with $A$ being the Hadamard matrix $H_\nu$, the ``no refinement'' versions of our policies would terminate according to the criterion $\|I_n-{1\over k}A_k^TA_k\|_\infty<{1\over 2s}$, which, on a closest inspection, is nothing but a slightly spoiled version of the goodness test based on mutual incoherence \cite{Donohoetal}\footnote{The mutual incoherence test is as follows: given a $k\times n$ matrix $B=[b_1,...,b_n]$ with nonzero columns, we compute the quantity $\mu(B)=\max\limits_{i\neq j}|b_i^Tb_j|/b_i^Tb_i$ and claim that $B$ is $s$-good for all $s$ such that $s<{1+\mu(B)\over 2\mu(B)}$. With the Hadamard $A$, the ``no refinement'' criterion for our scheme is nothing but $s<{1\over 2\mu(A^k)}$.}. In the experiments we are reporting, this criterion is essentially weaker that the one based on (\ref{relation}): for the best, over the 6 runs of the algorithms $\cA$ and $\cB$, 10-good submatrices $A_m$ of $H_{11}$ matrices we got the test based on mutual incoherence certifies the levels of goodness as low as 5 (in the case of $\cB$) and 7 (in the case of $\cA$).

\section{Low rank approximation of arbitrary matrices}\label{sec2}
\subsection{Randomized approximation}
\begin{proposition}\label{propinin}
Let $D\geq0$, and let $P=[p_1^T;...;p_m^T]\in\R^{m\times d}$ and $Q=[q_1^T;...;q_n^T]\in\R^{n\times d}$ be such that the Euclidean norms of the vectors $p_i$ and $q_j$  of $P$ and $Q$ are bounded by $\sqrt{D}$. Let an $m\times n$ matrix $A$ be represented as
$$
A=PQ^T
$$
\par
Given a positive integer $k$, consider the random matrix
\begin{equation}\label{eq1}
A_k={1\over k}P\left[\sum_{i=1}^k\xi_i\xi_i^T\right]Q^T={1\over k}\sum_{i=1}^k \eta_i\zeta_i^T,\,\,\eta_i:=\eta[\xi_i]=P\xi_i,\,\zeta_i:=\zeta[\xi_i]=Q\xi_i,
\end{equation}
where $\xi_i\sim \cN(0,I_d)$, $i=1,...,k$ are independent standard normal random vectors from $\R^d$. Then
\begin{equation}\label{then}
k\geq 8\ln(4mn)\Rightarrow \Prob\{\|A_k-A\|_\infty\leq {\sqrt{8\ln(4mn)}D\over \sqrt{k}}\}\geq {1\over 2}.
\end{equation}
\end{proposition}
For the proof, see Appendix \ref{appnew}.

\subsection{The norm associated with Proposition \ref{propinin}}\label{sectnorm}
Some remarks are in order. The result of Proposition \ref{propinin} brings to our attention to the smallest $D$ such that a given matrix $A$ can be decomposed into the product $PQ^T$ of two matrices with the Euclidean lengths of the rows not exceeding $\sqrt{D}$. On the closest inspection, $D$ turns out to be an easy-to-describe norm on the space $\R^{m\times n}$ of $m\times n$ matrices.
Specifically, let $\|A\|$, $A\in\R^{m\times n}$, be
$$
\|A\|=\min\limits_{t,M,N}\left\{t: \left[\begin{array}{c|c}M&A\cr\hline A^T&N\cr\end{array}\right]\succeq0,M_{ii}\leq t~\forall i, N_{jj}\leq t~\forall j\right\}
$$
This relation clearly defines a norm, and one clearly has $\|A\|=\|A^T\|$.
\begin{proposition}\label{prop1} For every $A\in\R^{m\times n}$, there exists representation $A=PQ^T$ with $P\in\R^{m\times (m+n)}$, $Q\in\R^{n\times (m+n)}$ and Euclidean norms of rows in $P,Q$ not exceeding $\sqrt{\|A\|}$.
Vice versa, if $A=PQ^T$ with the rows in $P,Q$ of Euclidean norms not exceeding $\sqrt{D}$, then $\|A\|\leq D$.
\end{proposition}
The next result summarizes the basic properties of the norm we have introduced.
\begin{proposition}\label{propsummary} Let $A$ be an $m\times n$ matrix. Then
\par
{\rm (i)}  $\|A\|_\infty\leq \|A\|\leq\sqrt{\min[m,n]}\|A\|_\infty$.\par
{\rm (ii)}  $\|A\|\leq \|A\|_{2,2}$, where $\|A\|_{2,2}$ is the usual spectral norm of $A$ (the maximal singular value).\par
{\rm (iii)}  If $A$ is symmetric positive semidefinite, then $\|A\|=\|A\|_\infty$.\par
{\rm (iv)} If the Euclidean norms of all rows (or all columns) of $A$ are $\leq D$, then $\|A\|\leq D$.
\end{proposition}
For the proof, see Appendix \ref{appnewnew}.
\subsection{Lower bound}
We have seen that if $A\in\R^{m\times n}$, then the $\|\cdot\|_\infty$-error of the best in this norm approximation of $A$ by a matrix of rank $k$ is {\sl at most} $O(1)\sqrt{\ln[mn]} \|A\|k^{-1/2}$.
We intend to demonstrate that in general this bound is unimprovable, up to a logarithmic in $m$ and $n$ factor. Specifically, the following result holds:
\begin{proposition}\label{proplb}
When $n\geq 2k$, the $\|\cdot\|_\infty$ error of any approximation of the unit matrix $I_n$ by a matrix of rank $k$ is at least
\begin{equation}\label{eq10}
{1\over 2\sqrt{k}}.
\end{equation}
Note that $\|I_n\|=1$.
\end{proposition}
{\bf Proof} [cf. \cite[Proposition 4.2]{CSFirst}] Let $\alpha(n,k)$ be the minimal $\|\cdot\|_\infty$ error of approximation of $I_n$ by a matrix of rank $\leq k$; this function clearly is nondecreasing in $n$. Let $\nu$ be an integer such that $k<\nu\leq n$, and $A$ be an $\nu\times\nu$ matrix of rank $\leq k$ such that $\|I_\nu-A\|_\infty=\alpha:=\alpha(\nu,k)$. By variational
characterization of singular values, at least $\nu-k$ singular values of $I_\nu-A$ are $\geq1$, whence $\Tr([I_\nu-A][I_\nu-A]^T)\geq\nu-k$. On the other hand, $\|I_\nu-A\|_\infty\leq\alpha$, whence $\Tr([I_\nu-A][I_\nu-A]^T)\leq\nu^2\alpha^2$. We conclude that $\alpha^2\geq{\nu-k\over \nu^2}$ for all $\nu$ with $k<\nu\leq n$, whence $\alpha^2\geq {1\over 4k}$ when $n\geq 2k$. \qed

We have seen that when $A\in\R^{m\times n}$, $A$ admits rank-$k$ approximations with the approximation error, measured in the $\|\cdot\|_\infty$-norm, of order of $\sqrt{\ln[mn]}\|A\|k^{-1/2}$. Note that the error bound deteriorates as $\|A\|$ grows. A natural question is, whether we could get similar results with a ``weaker''  norm of $A$ as a scaling factor. Seemingly the best we could hope for is $\|A\|_\infty$ in the role of the scaling factor, meaning that whenever all entries of an $m\times n$ matrix $A$ are in $[-1,1]$, $A$ can be approximated in $\|\cdot\|_\infty$-norm by a matrix of rank $k$ with approximation error which, up to a logarithmic in $m,n$ factor, depends solely on $k$ and goes to $0$ as $k$ goes to infinity. Unfortunately, the reality does not meet this hope. Specifically, let $A=H_\nu$ be the $n\times n$ Hadamard matrix ($n=2^\nu$), so that $\|A\|_\infty=1$.
Since $H^TH=n I_n$,  all $n$ singular values of the matrix are equal to $\sqrt{n}$, whence for every $n\times n$ matrix $B$  of rank $k<n$ the Frobenius norm of $A-B$ is at least $\sqrt{n(n-k)}$, meaning that the uniform norm of $A-B$ is at least $\sqrt{1-k/n}$. We conclude that the rank of a matrix which approximates $A$ with $\|\cdot\|_\infty$-error $\leq1/4$ should be of order of $n$.

\appendix
\section{Proof of Lemma \ref{dentropy}}\label{appa}
Properties (i) and (ii) are immediate consequences of the definition of $V_\beta$.
Observe that $V_\beta$ is convex and continuously differentiable with
\bse \left|{d\over dt}\big|_{t=0}V_\beta(x+th)\right|=\left|{\sum_{i=1}^d \sinh(x_i/\beta)h_i\over \sum_{i=1}^d
\cosh(x_i/\beta)}\right|\leq \|h\|_\infty\,\forall h, \ese
whence $\|V'_\beta(x)\|_1\le 1$ for $x\in \R^d$. Verification of \rf{2td} takes one line: $V_\beta$ is twice continuously differentiable with
\bse {d^2\over dt^2}\big|_{t=0}V_\beta(x+th)=\beta^{-1}{\sum_{i=1}^d \cosh(x_i/\beta)h_i^2\over \sum_{i=1}^d \cosh(x_i/\beta)}-\beta^{-1}{\left(\sum_{i=1}^d \sinh(x_i/\beta)h_i\right)^2\over
\left(\sum_{i=1}^d \cosh(x_i/\beta)\right)^2}\leq \beta^{-1}\|h\|_\infty^2. \ese \qed

\section{Problems (\ref{Opti}) in the case of Hadamard matrix $A$}\label{appHad}
We claim that if $A_m$ is an $m\times 2^\nu$ submatrix of the Hadamard matrix $H_\nu$ of order $n=2^\nu$, then the optimal values in all problems (\ref{Opti}) are equal to each other. The explanation is a s follows. Let $G$ be a finite abelian group of cardinality $n$. Recall that a {\sl character} of $G$ is a complex-valued function $\xi(g)$ such that $\xi(0)=1$ and $\xi(g+h)=\xi(g)\xi(h)$ for all $g,h\in G$; from this definition it immediately follows that $|\xi(g)|\equiv 1$. The characters of a finite abelian group $G$ form abelian group $G_*$, the multiplication being the pointwise multiplication of functions, and
this group is isomorphic to $G$. The Fourier Transform matrix associated with $G$ is the $n\times n$ matrix with rows indexed by $\xi\in G_*$, columns indexed by $g\in G$ and entries $\xi(g)$. For example, the usual DFT matrix of order $n$ corresponds to the cyclic group $G=\Z_n:=\Z/n\Z$, while the Hadamard matrix $H_\nu$ is nothing but the Fourier Transform matrix associated with $G=[\Z_2]^\nu$ (in this case, all characters take values $\pm1$). For $g\in G$ let $e_g(h)$ stands for the function on $G$ which is equal to 1 at $h=g$ and is equal to 0 at $h\neq g$. Given an $m$-element subset $Q$ of $G_*$, consider the submatrix
$A=[\xi(g)]_{{\xi\in Q\atop g\in G}}$ of the Fourier Transform matrix, along with $n$ optimization problems
$$
\min\limits_{y\in \C^m}\|\Re[e_g-A^Ty]\|_\infty=\min\limits_{y_\xi\in \C} \max\limits_{h\in G}|\Re[e_g(h)-\sum_{\xi\in Q} y_\xi\xi(h)]|\eqno{(P_g)}
$$
These problems clearly have equal optimal values, due to
$$
\begin{array}{l}
\max\limits_{h\in G}|\Re[e_g(h)-\sum_{\xi\in Q} y_\xi\xi(h)]| = \max\limits_{h\in G}|\Re[e_0(h-g)-\sum_{\xi\in Q} [y_\xi\xi(g)]\xi(h-g)]|\\ = \max\limits_{f=h-g\in G}|\Re[e_0(f)-\sum_{\xi\in
Q}[y_\xi\xi(g)]\xi(f)|.\\
\end{array}
$$
As applied to $G=\Z_2^\nu$, this observation implies that all quantities given by (\ref{Opti}) are the same.

\section{Proof of Proposition \ref{propinin}}\label{appnew}
The reasoning to follow is completely standard. Let us fix $i$, $1\leq i\leq m$, and $j$, $1\leq j\leq n$, and let $\xi\sim \cN(0,I_d)$,
$\mu=D^{-1/2}p^T_i\xi$, $\nu=D^{-1/2} q_{j}^T\xi$, and $\alpha=D^{-1}A_{ij}$. Then $[\mu;\nu]$ is a normal random vector with $\bE\{\mu^2\}\leq1$, $\bE\{\nu^2\}\leq1$ and $\bE\{\mu\nu\}=\alpha$. We can find a normal random vector $z=[u;v]\sim\cN(0,I_2)$ such that $\mu=au$, $\nu=bu+cv$; note that $a^2\leq1$, $b^2+c^2\leq1$ and $ab=\bE\{\mu\nu\}=\alpha$. Note that $\mu\nu=z^TBz$ with $B=\left[\begin{array}{cc}ab&ac/2\cr ac/2&0\cr\end{array}\right]$. Denoting $\lambda_1$, $\lambda_2$ the eigenvalues of $B$, we have
\begin{equation}\label{eq77}
\lambda_1+\lambda_2=\Tr(B)=ab=\alpha,~~~~\lambda_1^2+\lambda_2^2=\Tr(BB^T)=a^2(b^2+c^2/2)\leq 1.
\end{equation}
Now let $\gamma\in \R$ be such that $|\gamma|\leq1/4$. By (\ref{eq77}) we have $I_2-2B\succ0$, whence
$$
\bE\{\exp\{\gamma\mu\nu\}\}=\bE\{\exp\{\gamma z^TBz\}\}=\Det^{-1/2}(I_2-2\gamma B)=\left[(1-2\gamma\lambda_1)(1-2\gamma\lambda_2)\right]^{-1/2}.$$ Let $t=\sqrt{8\ln(4mn)}$ and $k\geq t$, and let
$[\mu_\ell;\nu_\ell]$, $1\leq\ell\leq k$, be independent random vectors with the same distribution as that of $[\mu;\nu]$. Then for every $\gamma\in(0,1/4]$ we have
\bse
\kappa_+
&:=&\Prob\{k[A_k]_{ij}>D[\alpha k+tk^{1/2}]\} =\Prob\{\sum_{\ell=1}^k\mu_\ell\nu_\ell\geq \alpha k+tk^{1/2}\}\\
&\leq& \bE\{\exp\{\gamma\sum_{\ell=1}^k\mu_\ell\nu_\ell\}\}\exp\{-\gamma k(\alpha+k^{-1/2}t)\} \\ &=&\left[\bE\{\exp\{\gamma\mu\nu\}\}\right]^k\exp\{-\gamma k(\alpha + k^{-1/2}t)\}\\ &=&\left[(1-2\gamma\lambda_1)(1-2\gamma\lambda_2)\right]^{-k/2}\exp\{-\gamma k(\alpha+k^{-1/2}t)\},
\ese
so that
\bse
\ln\kappa_+ &\leq& {k\over 2}\left[-2\gamma(\alpha+k^{-1/2}t)-\ln(1-2\gamma\lambda_1)-\ln(1-2\gamma\lambda_2)\right]\\ &\leq& {k\over 2}\left[-2\gamma[\lambda_1+\lambda_2]-2\gamma k^{-1/2}t-\ln(1-2\gamma\lambda_1)-\ln(1-2\gamma\lambda_2)\right]\\
&\leq&{k\over 2}\left[-2\gamma k^{-1/2}t+4\gamma^2(\lambda_1^2+\lambda_2^2)\right]
\ese
where the last inequality follows from $|2\gamma\lambda_s|\leq 1/2$, for $s=1,2$, and $-\ln(1-r)-r\leq r^2$ when $|r|\leq1/2$.
Using \eqref{eq77} we obtain,
$$
\ln\kappa_+\leq {k\over 2}\left[-2\gamma k^{-1/2}t+4\gamma^2\right].
$$
Setting $\gamma={t\over 4k^{1/2}}$ (this results in $0<\gamma\leq1/4$ due to $k^{1/2}\geq t$), we get
$$
\Prob\{k[A_k]_{ij}> A_{ij} k+Dt{ k^{1/2}}\}=\kappa_+\leq \exp\{-t^2/8\}=(4mn)^{-1}.
$$
Letting $\kappa_-=\Prob\{k[A_k]_{ij}<A_{ij} k - Dk^{1/2}t\}$, we have
$$\kappa_-\leq \bE\{\exp\{-\gamma \sum_{\ell=1}^k\mu_\ell\nu_\ell\}\}\exp\{-\gamma k(-\alpha+{ k^{-1/2}}t)\}$$
for all $\gamma\in(0,1/4]$, whence, same as above,
$$\Prob\{k[A_k]_{ij}< kA_{ij}-Dk^{1/2}t\}=\kappa_-\leq (4mn)^{-1}.$$
We see that
$$
\Prob\{|[A_k]_{ij}-A_{ij}|>Dtk^{-1/2}\}\leq{1\over 2mn}.
$$
Since this relation holds  true for all $i,j$, we conclude that
$$
\Prob\{\|A_k-A\|_\infty> Dk^{-1/2}t\}\leq 1/2. \eqno{\hbox{\qed}}
$$
\section{Proofs for section \ref{sectnorm}}\label{appnewnew}
\paragraph{Proof of Proposition \ref{prop1}.}
First claim: there exist $M,N$ such that the matrix $\cA=\left[\begin{array}{c|c}M&A\cr\hline A^T&N\cr\end{array}\right]$ is positive semidefinite and has all diagonal entries, and then all entries, in $[-\|A\|,\|A\|]$. Let $\cA=\cB\cB^T$; then the rows in $\cB$ have Euclidean norms $\leq \sqrt{\|A\|}$. Representing $\cB=[P;Q]$ with $m$ rows in $P$ and $n$ rows in $Q$, the relation $[P;Q][P;Q]^T=\cA$ implies that $A=PQ^T$.
\par Second claim: Let $A=PQ^T$ with the Euclidean norms of rows in $P,Q$ not exceeding $\sqrt{D}$. Then $0\preceq \left[\begin{array}{c}P\cr
Q\cr\end{array}\right]\left[\begin{array}{c}P\cr Q\cr\end{array}\right]^T=\left[\begin{array}{c|c}PP^T&A\cr\hline A^T&QQ^T\cr\end{array}\right]$ and the diagonal entries in $M=PP^T$ and $N=QQ^T$ do not
exceed $D$. \qed
\paragraph{Proof of Proposition \ref{propsummary}.} {\bf (i)}:  The first inequality in (i) is evident. Let us prove the second. W.l.o.g. we can assume $\|A\|_\infty\leq1$. In this case our statement reads
$$
\Opt:=\min\limits_{t,M,N}\left\{t: \left[\begin{array}{c|c}M&A\cr\hline A^T&N\cr\end{array}\right]\succeq0,t-M_{ii}\geq0\forall i, t-N_{jj}\geq0\,\forall j\right\}\leq D=\sqrt{\min[m,n]}.
$$
Assume, on the contrary, that $\Opt>D$. Since the semidefinite problem defining $\Opt$ is strictly feasible, the dual problem
$$
\max\limits_{X,Y,Z,\lambda,\rho}\left\{-2\Tr(Z^TA):\begin{array}{l} \left[\begin{array}{c|c}X&Z\cr\hline Z^T&Y\cr\end{array}\right]\succeq0\\ \lambda\geq0,\rho\geq0,\sum_i\lambda_i+\sum_j\rho_j=1\\
\Tr(XM)+\Tr(YN)+\sum_i\lambda_i(t-M_{ii})\\ +\sum_j\rho_j(t-N_{jj})\equiv t\,\forall M,N,t\\
\end{array}\right\}
$$
has a feasible solution with value of the objective $>D$. In other words, there exist nonnegative vectors $\lambda\in\R^m$, $\rho\in\R^n$ and a matrix $V=-Z\in\R^{m\times n}$ such that
$$
\begin{array}{ll}
(a)&\left[\begin{array}{c|c}\Diag\{\lambda\}&V\cr\hline V^T&\Diag\{\rho\}\cr\end{array}\right]\succeq0\\ (b)&\sum_i\lambda_i+\sum_j\rho_j=1\\ (c)&2\Tr(V^TA)>D.\\
\end{array}
$$
By $(a)$, letting $L=\Diag\{\sqrt{\lambda_i}\}$, $R=\Diag\{\sqrt{\rho_j}\}$, we have $V=LWR$ with certain $W$, $\|W\|_{2,2}\leq1$ ($\|\cdot\|_{2,2}$ is the usual matrix norm, the maximum singular value), thus
\begin{equation}\label{eq101}
\begin{array}{l}
2\Tr(V^TA)=2\Tr(RW^TLA)\leq 2\sum_{i,j}|[RW^TL]_{ij}|=2\sum\limits_{i,j}L_{ii}|W_{ij}|R_{jj}\\ =2\sum_iL_{ii}\sum_{j}|W_{ij}|R_{jj}\leq2\|[|W|]_{i,j}\|_{2,2}\sqrt{\sum_iL_{ii}^2}\sqrt{\sum_jR_{jj}^2}\\
\underbrace{\leq}_{(*)} 2\sqrt{\min[m,n]}\sqrt{(\sum_i\lambda_i)(\sum_j\rho_j)}\leq D,\\
\end{array}
\end{equation}
where the concluding $\leq$ is due to $(b)$, and $(*)$ is given by the following reasoning: w.l.o.g. we can assume that $n\leq m$. Since $W$ is of the matrix norm $\leq1$, the columns $U_j$ of $U=[|W_{ij}|]_{i,j}$ satisfy $\|U_j\|_2\leq1$, whence
$$
\|Ux\|_2\leq \sum_{i=1}^n |x_j|\|U_j\|_2\leq\sqrt{n}\|x\|_2\,\,\forall x.
$$
The resulting inequality in (\ref{eq101}) contradicts $(c)$; we have arrived at a desired contradiction. (i) is proved.\\
{\bf (ii)}: This is evident, since $\left[\begin{array}{c|c}\|A\|_{2,2}I_m&A\cr\hline A^T&\|A\|_{2,2}I_n\cr\end{array}\right]\succeq0$.\\
{\bf (iii)}: This is evident, since for $A\succeq0$ we have $\left[\begin{array}{c|c}A&A\cr\hline A&A\cr\end{array}\right]\succeq0$. \\
{\bf (iv)}:  Since $\|A\|=\|A^T\|$, it suffices to consider the case when the rows of $A$ are of the norm not exceeding $D$. In this case, the result is readily given by the fact that
$\left[\begin{array}{c|c}D^{-1}AA^T&A\cr\hline A^T&DI_n\cr\end{array}\right]\succeq0$. \qed
\end{document}